\documentclass[11pt]{article}

\usepackage{amsmath}
\usepackage{amsthm}
\usepackage{latexsym}
\usepackage{amssymb}
\usepackage{graphicx}
\usepackage{amsmath}
\usepackage{amssymb}
\usepackage{graphicx}
\usepackage{amsmath}
\usepackage{amsfonts}
\usepackage{euscript}
\usepackage{epsfig}
 \usepackage{color}
\usepackage{amsmath}
\usepackage{latexsym}
\usepackage{amssymb}
\usepackage{cite}
\usepackage{indentfirst}
\usepackage{enumerate}
\usepackage{tikz}

\newtheorem{lem}{Lemma}[section]

\newtheorem{thm}{Theorem}[section]
\newtheorem{cor}{Corollary}[section]

\title{\bf Generalized Weyl Conformal Curvature Tensor of
Generalized Riemannian Space}
\author{Nenad O. Vesi\'c}
\date{}

\numberwithin{equation}{section}

\makeatletter
  \def\my@tag@font{\normalsize}
  \def\maketag@@@#1{\hbox{\m@th\normalfont\my@tag@font#1}}
  \let\amsmath@eqref\eqref
  \renewcommand\eqref[1]{{\let\my@tag@font\relax\amsmath@eqref{#1}}}
\makeatother

\begin{document}
  \maketitle

  \begin{abstract}
    It is generalized Weyl
    conformal curvature tensor in the case of a conformal mappings
    of a generalized Riemannian space in this paper. Moreover, it is found
    universal generalizations of it without any additional
    assumption. A method used in this paper may help different
    scientists in their researching.\\[3pt]


    \noindent\textbf{$2010$ Math Subj. Classification:} 14L24, 14A22,
    53A30, 53B20
  \end{abstract}

  \section{Introduction}

{Many research papers, books and monographs are dedicated to
  development of the theory of conformal mappings and its
  applications.} Some of authors who
  have contributed to this development are H. M. Abood \cite{wc1}, S.
  Bochner \cite{wc2}, L. P. Eisenhart \cite{wc4}, \linebreak S. B. Mathur \cite{wc7}, Josef Mike\v s with his
  research group \cite{mikesc1, wc3, mikes2, wc15, mikes1, zlatconform, marijaconform,
  mikesc2}, \linebreak S. M.
  Min\v ci\'c \cite{wc14}, P. Mocanu
  \cite{wc11}, M. Prvanovi\'c \cite{wc12}, N. S. Sinyukov \cite{wc13},
  M. Lj. Zlatanovi\'c, M. Najdanovi\'c \cite{zlatconform, marijaconform} and many others.
   A. Einstein \cite{e1, e2,
  e3} based the theory of general relativity
  on non-symmetric affine connection.
 E. Goulart and M. Novello \cite{wc5} such as H. Zhang, Y.
  Zhang, X-Z. Li \cite{zhang1} applied the theory of conformal mappings in
  physics.

  The main purpose of this paper is to make analogies between
  invariants of geodesic and conformal mappings, i.e. we want to
  examine are there analogies of Thomas projective parameter and
  Weyl projective tensor as invariants of conformal mappings in
  here.

  \subsection{Generalized Riemannian spaces}

  %An $N$-dimensional manifold $\mathcal M_N$ endowed with metric
  %tensor $g_{ij}$ symmetric by indices $i$ and $j$ is the Riemannian space $\mathbb R_N$ \cite{mikesc1, wc3, mikes2, wc15, mikes1,
  %mikesc2}.
  {Based on the Eisenhart's results \cite{wc4}, many
  authors started the researches about conformal mappings between
  Riemannian and
  generalized Riemannian spaces as well as about their invariants
   (see \cite{wc1, wc2, wc5, mikesc2, wc3, mikesc1, mikes2, wc15, mikes1, wc9, wc10,
   marijaconform, wc14, FilomatConform, zlatconform,
  wc12})}.\vspace{.2cm}

   An $N$-dimensional manifold $\mathcal M_N$ endowed with a
  metric tensor $G_{ij}$ non-symmetric by indices $i$ and $j$, is the
  generalized Riemannian space $\mathbb{GR}_N$ {\cite{wc4}}.
  The symmetric and anti-symmetric
  parts of the metric tensor $G_{ij}$ are respectively defined as

      \begin{eqnarray}
      g_{{ij}}=\frac12(G_{ij}+G_{ji})&\mbox{and}&
      F_{{ij}}=\frac12(G_{ij}-G_{ji}).
      \label{eq:gsymasym}
    \end{eqnarray}

\pagebreak
  \noindent (Generalized) Christoffel symbols of the space $\mathbb{GR}_N$ are

    \begin{eqnarray}
    \Gamma_{i.jk}=\frac12(G_{ji,k}-G_{jk,i}+G_{ik,j})&\mbox{and}&
    \Gamma^i_{jk}=g^{{i\alpha}}\Gamma_{\alpha.jk},
  \end{eqnarray}

  \noindent for partial derivative $\partial/\partial x^k$ denoted
  by comma. Christoffel symbols $\Gamma^i_{jk}$ are a case of non-symmetric linear connection, where their
  symmetric and anti-symmetric parts are given by

  \begin{align}
    &\gamma{}^i_{{jk}}=\frac12(\Gamma^i_{jk}+
    \Gamma^i_{kj})=\frac12g^{{i\alpha}}(g_{{j\alpha},k}-
    g_{{jk},\alpha}+g_{{\alpha k},j}),\label{eq:gammasym}\\&
    T^i_{{jk}}=\frac12(\Gamma^i_{jk}-
    \Gamma^i_{kj})=\frac12g^{{i\alpha}}(F_{{j\alpha},k}-
    F_{{jk},\alpha}+F_{{\alpha
    k},j}).\label{eq:gammaasym}
  \end{align}

  \noindent The anti-symmetric part $T_{{jk}}^i$ is
  torsion tensor of the space $\mathbb{GR}_N$. It also
  holds

  \begin{eqnarray}
  \gamma{}^\alpha_{{j\alpha}}=\frac12\big(\ln{|g|}\big)_{,j}&
  \mbox{and}&T^\alpha_{j\alpha}=0,
  \end{eqnarray}

  \noindent for $g=\det[g_{{ij}}]\neq0$.
  The Riemannian space $\mathbb R_N$ endowed with the affine
  connection coefficients
  $\gamma{}^i_{{jk}}$ is called \emph{the associated space of
  $\mathbb{GR}_N$} {\cite{wc9, wc10,
  marijaconform, wc14, zlatconform}}.

Since the associated space $\mathbb R_N$ is usual Riemannian space,
there exists only one kind of covariant differentiation with regard
to the affine connection of this space:

    \begin{equation}
    a^i_{j;k}=a^i_{j,k}+\gamma{}^i_{{\alpha k}}a^\alpha_j-
    \gamma{}^\alpha_{{jk}}a^i_\alpha,\label{eq:covdev}
  \end{equation}

  \noindent where  $a^i_j$ is a tensor of the type $(1,1)$.
  It also exists only one curvature tensor

{
  \begin{equation}
  R^i_{jmn}=\gamma{}^i_{{jm},n}-
  \gamma{}^i_{{jn},m}+
  \gamma{}^\alpha_{{jm}}\gamma{}^i_{{\alpha n}}-
  \gamma{}^\alpha_{{jn}}\gamma{}^i_{{\alpha m}},
  \label{eq:r}
\end{equation}
}

\noindent of the associated space $\mathbb R_N$.

In the generalized Riemannian space $\mathbb{GR}_N$, one can
consider four kinds of covariant differentiation \cite{wc9, wc10}

\begin{eqnarray}
  a^i_{j\underset1|k}=a^i_{j,k}+\Gamma^i_{\alpha k}a^\alpha_j-
  \Gamma^\alpha_{jk}a^i_\alpha&&
  a^i_{j\underset2|k}=a^i_{j,k}+\Gamma^i_{k\alpha}a^\alpha_j-
  \Gamma^\alpha_{kj}a^i_\alpha,\label{eq:covdev1}\\
  a^i_{j\underset3|k}=a^i_{j,k}+\Gamma^i_{\alpha k}a^\alpha_j-
  \Gamma^\alpha_{kj}a^i_\alpha&&
  a^i_{j\underset4|k}=a^i_{j,k}+\Gamma^i_{k\alpha}a^\alpha_j-
  \Gamma^\alpha_{jk}a^i_\alpha.\label{eq:covdev4}
\end{eqnarray}

\noindent Also, in the generalized Riemannian space $\mathbb{GR}_N$,
there exist twelve curvature tensors {\cite{wc9}}. These curvature
tensors are elements of the family

\begin{equation}
    K^i_{jmn}=R^i_{jmn}+uT^i_{{jm};n}+
  u'T^i_{{jn};m}+
  vT^\alpha_{{jm}}T^i_{{\alpha
  n}}+v'T^\alpha_{{jn}}T^i_{{\alpha
  m}}+wT^\alpha_{{mn}}T^i_{{\alpha
  j}},\label{eq:R}
\end{equation}

\noindent for the corresponding $u,u',v,v',w\in\mathbb R$. It is
proved that five of these curvature tensors are linearly independent
{\cite{wc10}}.

\subsection{Conformal mappings of generalized Riemannian space}

A mapping $f:\mathbb{GR}_N\rightarrow\mathbb{G\overline R}_N$
determined by the equation

\begin{equation}
  \overline{G}_{ij}=e^{2\psi}G_{ij},
  \label{eq:conformGRNdfn}
\end{equation}

\noindent for a scalar function $\psi$ is the conformal mapping of
 space $\mathbb{GR}_N$ \cite{wc14}. The basic equation of this mapping is

 \begin{equation}
   \overline\Gamma^i_{jk}=
   \Gamma^i_{jk}+\psi_j\delta^i_k+\psi_k\delta^i_j-\psi^ig_{{jk}}+\xi^i_{jk},
   \label{eq:conformGRNbasic}
 \end{equation}

 \noindent for $\psi_i=\partial\psi/\partial x^i,\psi^i=g^{i\alpha}\psi_\alpha$,
 and tensor
 $\xi^i_{jk}$ anti-symmetric by indices $j$ and $k$. After
 antisymmetrization of the basic equation (\ref{eq:conformGRNbasic})
 by indices $j$ and $k$, we obtain that is

 \begin{equation*}
   \xi^i_{jk}=\overline T^i_{{jk}}-
   T^i_{{jk}}.
 \end{equation*}

 For this reason, it were studied equitorsion conformal mappings
 (the  case of $\xi^i_{jk}=0$)\linebreak \cite{zlatconform,
wc14, marijaconform}. In last author's research
\cite{FilomatConform}, it were studied conformal mappings which do
not preserve torsion tensor.

\subsection{Motivation}

{Geometrical objects that are invariant with respect to conformal
mappings play an important role in the theory of gravity
\cite{zhang1, wc5}}.

Let $f:\mathbb{GR}_N\rightarrow\mathbb{G\overline R}_N$ be a
conformal mapping of generalized Riemannian space $\mathbb{GR}_N$.
Weyl conformal curvature tensor

{
\begin{align}
  &\aligned
  C^i_{jmn}&=R^i_{jmn}+\frac1{N-2}(\delta^i_nR_{jm}-
      \delta^i_mR_{jn})+\frac1{N-2}(g_{{jm}}
      R^i_n-g_{{jn}}R^i_m)\\&+
      \frac{
      R}{(N-1)(N-2)}(\delta^i_mg_{{jn}}-\delta^i_ng_{{jm}}),
  \endaligned\label{eq:Weylconform}\\&
  \aligned
  C_{ijmn}&=R_{ijmn}+\frac1{N-2}(g_{in}R_{jm}-
      g_{im}R_{jn})+\frac1{N-2}(g_{{jm}}
      R_{in}-g_{{jn}}R_{im})\\&+
      \frac{
      R}{(N-1)(N-2)}(g_{im}g_{{jn}}-g_{in}g_{{jm}}),
  \endaligned\tag{1.13'}\label{eq:cWeylcpmfprm}
\end{align}
}

\noindent are invariants of this mapping obtained from the change of
curvature tensor $R^i_{jmn}$.

S. Bochner (see \cite{wc2}) generalized the covariant Weyl conformal
curvature tensor. This generalization is

\begin{equation}
  \aligned
  B_{ij^\ast mn^\ast}&=K_{ij^\ast
  mn^\ast}\\&-\frac1{k+2}\big(G_{ij^\ast}K_{mn^\ast}+
  G_{in^\ast}K_{mj^\ast}+G_{mj^\ast}K_{in^\ast}+
  G_{mn^\ast}K_{ij^\ast}\big)\\&+\frac K{2(k+1)(k+2)}\big(
  G_{ij^\ast}G_{mn^\ast}+G_{in^\ast}G_{j^\ast m}\big),
  \endaligned
  \label{eq:bochnert}
\end{equation}

\noindent from the change of covariant curvature tensor

\begin{equation*}
  K_{ij^\ast mn^\ast }=T_{i.j^\ast m;n^\ast}-
  T_{i.j^\ast n^\ast;m}+T^\alpha_{j^\ast m}T_{i.\alpha
  n^\ast}-T^\alpha_{j^\ast n^\ast}T_{i.\alpha m},
\end{equation*}

\noindent of $N=2k$-dimensional space with Hermitian metrics

\begin{equation*}
  ds^2=2g_{\alpha\beta^\ast}dz_\alpha d\overline z_\beta
\end{equation*}

\noindent which satisfies the K\"ahlerian assumption

\begin{equation*}
  \frac{\partial g_{\alpha\gamma^\ast}}{\partial z_\beta}=
  \frac{\partial g_{\beta\gamma^\ast}}{\partial z_\alpha}.
\end{equation*}

\section*{Purposes of the paper}

 {The family of invariants of a
conformal mapping which preserves the torsion tensor
 obtained from the change of the family of curvature tensors
(\ref{eq:R}) is searched in \cite{wc14}}.

  {Recently, we obtained that
  it is not necessary to assume the equitorsioness of
  conformal mappings to find their invariants \cite{FilomatConform}.
  In that article, it is obtained one family of generalizations of
  Weyl conformal curvature tensor. In this manuscript, we wish to find some other families of invariants of
  conformal mappings of a generalized Riemannian space.}

  The aims of this paper are:

  \begin{enumerate}
    \item To obtain necessary and sufficient conditions for a
    mapping $f:\mathbb{GR}_N\rightarrow\mathbb{G\overline R}_N$ to
    be a conformal one.
    \item To find invariants of a conformal mapping
    $f:\mathbb{GR}_N\rightarrow\mathbb{G\overline R}_N$ obtained
    from changes of Christoffel symbols $\Gamma^i_{jk}$ under this
    mapping.
    \item To generalize Weyl conformal curvature tensors (\ref{eq:Weylconform},
    \ref{eq:cWeylcpmfprm}).
   \end{enumerate}

  \section{Main results}

  Before the main examinations, we may notice that invariants of geometric
  mappings have been obtained from changes of affine connection
  coefficients or from  changes of curvature tensors of the space $\mathbb{GR}_N$
   under these
  mappings. Invariants obtained from changes of affine connection
  coefficients of the space $\mathbb{GR}_N$
  under a mapping will be called \textit{the conformal invariants of Thomas type}.
  Invariants obtained from changes of
  curvature tensors of the space $\mathbb{GR}_N$ under a mapping
   will be called \textit{the conformal
  invariants of Weyl type}.

%\pagebreak

  \subsection{Conformal invariants of Thomas type}

  Let $f:\mathbb{GR}_N\rightarrow\mathbb{G\overline R}_N$ be
  a conformal mapping of the space $\mathbb{GR}_N$. After symmetrize
  the equation
  (\ref{eq:conformGRNbasic}) by indices $j$ and $k$ and contract
  the symmetrized equation by $i$ and $k$, we established the correctness
  of the following equalities:

  \begin{equation}
    \psi_j=\frac1N({\overline\gamma}{}^\alpha_{{j\alpha}}-
    \gamma{}^\alpha_{{j\alpha}})=\frac1{2N}\Big(
    \big(\ln{|\overline g|}\big)_{,j}-\big(\ln{|g|}\big)_{,j}\Big).
    \label{eq:conformpsi}
  \end{equation}

  \noindent In this way, we proved that the difference
  $p{}^i_{jk}={\overline\gamma}{}^i_{jk}-
  \gamma{}^i_{jk}$ is

  \begin{equation}
    p{}^i_{{jk}}={\overline\zeta}^i_{(1)jk}-
    \zeta^i_{(1)jk}={\overline\zeta}^i_{(2)jk}-
    \zeta{}^i_{(2)jk},\label{eq:Pconform12}
  \end{equation}

  \noindent for
  \begin{eqnarray*}
    \zeta^i_{(1)jk}=\gamma^i_{{jk}}&\mbox{and}&
    \zeta^i_{(2)jk}=\frac1{2N}\Big(\big(\ln{|g|}\big)_{,j}\delta^i_k+
    \big(\ln{|g|}\big)_{,k}\delta^i_j-\big(\ln{|g|}\big)_{,\alpha}g^{{i\alpha}}
    g_{{jk}}
    \Big)
  \end{eqnarray*}

  \noindent and the corresponding
  ${\overline\zeta}^i_{(1)jk}$ and
  ${\overline\zeta}^i_{(2)jk}$.

  Because it holds $\overline
  g^{{ij}}=e^{-2\psi}g^{{ij}}$, we have that is

    \begin{equation}
    \overline g^{{ij}}\overline G_{mn}=
    e^{-2\psi}g^{{ij}}e^{2\psi}G_{mn}=g^{{ij}}G_{mn},
    \label{eq:ggINV}
  \end{equation}

  \noindent i.e. $g^{{ij}}G_{mn}$ is an invariant of the
  mapping $f$.

  The torsion tensor $T^i_{jk}$ may be expressed as

    \begin{equation}
  \aligned
    T^i_{{jk}}&=\frac12g^{{i\alpha}}
    \big(F_{{j\alpha};k}-
    F_{{jk};\alpha}+
    F_{{\alpha k};j}\big)
    =\frac12\Big(\big(g^{{i\alpha}}
    F_{{j\alpha}}\big)_{;k}-
    \big(g^{{i\alpha}}F_{{jk}}\big)_{;\alpha}+
    \big(g^{{i\alpha}}F_{{\alpha
    k}}\big)_{;j}\Big).
  \endaligned\label{eq:Torsion}
  \end{equation}

  \noindent From this expression and the definition
  of the covariant derivation with respect to linear
  connection of the associated space $\mathbb R_N$
  given by the equation (\ref{eq:covdev}), we directly obtain that is

  \begin{equation}
    \aligned
    \overline T^i_{{jk}}-
    T^i_{{jk}}&=\frac12\big({\overline\gamma}{}^i_{{\alpha
    k}}\overline g^{{\alpha\beta}}\overline
    F_{{j\beta}}-{\overline\gamma}{}^i_{{\beta\alpha}}
    \overline g^{{\beta\alpha}}\overline
    F_{{jk}}+
    {\overline\gamma}{}^i_{{\alpha j}}\overline
    g^{{\beta\alpha}}\overline F_{{\beta k}}
    +
    {\overline\gamma}{}^\alpha_{{j\beta}}\overline
    g^{{i\beta}}\overline F_{{\alpha k}}+
    {\overline\gamma}{}^\alpha_{{k\beta}}\overline
    g^{{i\beta}}\overline F_{{j\alpha}}
    \big)\\&-
    \frac12\big(\gamma{}^i_{{\alpha
    k}}g^{{\alpha\beta}}
    F_{{j\beta}}-\gamma{}^i_{{\beta\alpha}}
    g^{{\beta\alpha}}
    F_{{jk}}+
    \gamma{}^i_{{\alpha j}}
    g^{{\beta\alpha}}F_{{\beta k}}
    +
    \gamma{}^\alpha_{{j\beta}}
    g^{{i\beta}}F_{{\alpha k}}+
    \gamma{}^\alpha_{{k\beta}}
    g^{{i\beta}}F_{{j\alpha}}
    \big)\\&\overset{(\ref{eq:ggINV})}=
    \frac12\big(p{}^i_{{\alpha
    k}}g^{{\alpha\beta}}
    F_{{j\beta}}-p{}^i_{{\beta\alpha}}
    g^{{\beta\alpha}}
    F_{{jk}}+
    p{}^i_{{\alpha j}}
    g^{{\beta\alpha}}F_{{\beta k}}
    +
    p{}^\alpha_{{j\beta}}
    g^{{i\beta}}F_{{\alpha k}}+
    p{}^\alpha_{{k\beta}}
    g^{{i\beta}}F_{{j\alpha}}
    \big).
    \endaligned\label{eq:Torsiontransformations}
  \end{equation}

  Based on the equations
  (\ref{eq:conformpsi}, \ref{eq:Pconform12}, \ref{eq:Torsiontransformations}) and the invariance
  (\ref{eq:ggINV}), we obtain that torsion tensors
  $T^i_{{jk}}$ and $\overline T^i_{{jk}}$
  satisfy the following relations:

  \begin{equation}
    \overline T^i_{{jk}}=
    T^i_{{jk}}+\tau^i_{(r)jk}-
   {\overline\tau}^i_{(r)jk},
    \label{eq:TtorsionC}
  \end{equation}

  \noindent for $r=(r_1,\ldots,r_5)\in\{1,2\}^5$ and

\begin{footnotesize}
  \begin{align}
    &\aligned
    {\overline\tau}^i_{(r)jk}&=-
    \frac12\big({\overline\zeta}^i_{(r_1)\alpha k}
    \overline g^{{\alpha\beta}}\overline
    F_{{j\beta}}-
    {\overline\zeta}^i_{(r_2)\beta\alpha}
    \overline g^{{\beta\alpha}}\overline
    F_{jk}+{\overline\zeta}^i_{(r_3)\alpha
    j}\overline g^{{\beta\alpha}}\overline
    F_{{\beta
    k}}+{\overline\zeta}^\alpha_{(r_4)j\beta}\overline
    g^{{i\beta}}\overline F_{{\alpha
    k}}+{\overline\zeta}{}^\alpha_{(r_5)k\beta}\overline
    g^{{i\beta}}
    \overline F_{{j\alpha}}
    \big),
    \endaligned\label{eq:ttauc}\\&
    \aligned
    \tau^i_{(r)jk}=-\frac12\big({\zeta}{}^i_{(r_1)\alpha k}
    g^{{\alpha\beta}}
    F_{{j\beta}}-
    {\zeta}{}^i_{(r_2)\beta\alpha}
    g^{{\beta\alpha}}
    F_{{jk}}+{\zeta}{}^i_{(r_3)\alpha
    j}g^{{\beta\alpha}}
    F_{{\beta
    k}}+{\zeta}{}^\alpha_{(r_4)j\beta}
    g^{{i\beta}}F_{{\alpha
    k}}+{\zeta}{}^\alpha_{(r_5)k\beta}
    g^{{i\beta}}
    F_{{j\alpha}}
    \big).
    \endaligned\label{eq:tauc}
  \end{align}
\end{footnotesize}

  The equations (\ref{eq:conformGRNbasic}, \ref{eq:conformpsi}, \ref{eq:TtorsionC}) prove that is

  \begin{eqnarray*}
    {{\overline{\mathcal T}}}{}^i_{(r)jk}=
    {{\mathcal T}}{}^i_{(r)jk},
  \end{eqnarray*}

  \noindent for
  \begin{equation}
  \aligned
    {\mathcal
    T}{}^i_{(r)jk}&=\Gamma^i_{jk}-\frac1{2N}\Big(\big(\ln{|g|}\big)_{,j}\delta^i_k+
    \big(\ln{|g|}\big)_{,k}\delta^i_j-\big(\ln{|g|}\big)_{,\alpha}
    g^{{i\alpha}}g_{{jk}}
    \Big)+\tau{}^i_{(r)jk},
  \endaligned\label{eq:Thomas}
  \end{equation}

  \noindent a random $r\in\{1,2\}^5$, such as for the corresponding ${\overline{\mathcal
  T}}{}^i_{(r)jk}$.

\pagebreak
  It holds the
  following lemma:

  \begin{lem}
    Let $f:\mathbb{GR}_N\rightarrow\mathbb{G\overline R}_N$ be a
    mapping of generalized Riemannian space $\mathbb{GR}_N$. The
    following statements are equivalent:

    \begin{enumerate}[\emph{1.}]
     \item[\emph{1.}] Basic equation of the mapping $f$ is
      \begin{equation}
        \aligned
        \overline\Gamma^i_{jk}&=
        \Gamma^i_{jk}+\frac1{2N}\Big(\big(\ln{|\overline g|}\big)_{,j}
      \delta^i_k+\big(\ln{|\overline g|}\big)_{,k}\delta^i_j-
      \big(\ln{|\overline g|}\big)_{,\alpha}\overline g^{{i\alpha}}
      \overline g_{{jk}}\Big)\\&-
      \frac1{2N}\Big(\big(\ln{|g|}\big)_{,j}
      \delta^i_k+\big(\ln{|g|}\big)_{,k}\delta^i_j-
      \big(\ln{|g|}\big)_{,\alpha}g^{{i\alpha}}
      g_{{jk}}\Big)+\tau{}^i_{(r)jk}-
      {\overline\tau}{}^i_{(r)jk},
        \endaligned\label{eq:conformbasic1}
      \end{equation}

      for some $r\in\{1,2\}^5$ and the
      corresponding $\tau{}^i_{(r)jk},
      {\overline\tau}{}^i_{(r)jk}$ defined by the
      equations \emph{(\ref{eq:ttauc}, \ref{eq:tauc})}.
    \item[\emph{2.}] For $r\in\{1,2\}^5$, the geometrical object
    ${\mathcal T}{}^i_{(r)jk}$, defined by the equation
    \emph{(\ref{eq:Thomas})} is an invariant of the mapping $f$.

    \item[\emph{3.}] Mapping $f$ is a conformal mapping.\qed
    \end{enumerate}
  \end{lem}

  \noindent The invariants
  ${\mathcal T}{}^i_{(r)jk}$, defined by the equation (\ref{eq:Thomas}), are \textbf{the $r$-th
  conformal invariants of Thomas type} of the
  conformal mapping $f$.

  \subsection{Invariants of Weyl type}

  From the family (\ref{eq:R}) of curvature tensors
  of the space $\mathbb{GR}_N$, we get that is

  \begin{align}
    &\aligned
    R^i_{jmn}=K^i_{jmn}-uT^i_{{jm};n}-
    u'T^i_{{jn};m}-
    vT^\alpha_{{jm}}T^i_{{\alpha
    n}}-
    v'T^\alpha_{{jn}}T^i_{{\alpha
    m}}-
    wT^\alpha_{{mn}}T^i_{{\alpha
    j}},
    \endaligned\label{eq:Rtenzor=Ktenzor}\\
    &\aligned
    R_{ij}=K_{ij}-uT^\alpha_{{ij};\alpha}-
    (v'+w)T^\alpha_{{i\beta}}
    T^\beta_{{\alpha j}},
    \endaligned\label{eq:Ricc=KRicc}\\
    &\aligned
    R^i_j=K^i_j-ug^{{i\alpha}}T^\beta_{{\alpha
    j};\beta}-(v'+w)g^{{i\alpha}}T^\beta_{{\alpha\gamma}}
    T^\gamma_{{\beta j}},
    \endaligned\label{eq:Rij=Kij}\\
    &\aligned
    R=K-(v'+w)T^\alpha_{{\gamma\beta}}
    T^\beta_{{\alpha\delta}}g^{{\gamma\delta}}.
    \endaligned\label{eq:R=K}
  \end{align}

\noindent  Based on these expressions and because Weyl conformal
curvature tensor (\ref{eq:Weylconform}) is an invariant of the
conformal mapping $f$, we obtain that it holds

  \begin{equation}
    \aligned
    {\overline R}{}^i_{jmn}&=R^i_{jmn}+
    \frac1{N-2}\big(\delta^i_mK_{jn}-\delta^i_n
    K_{jm}+ K^i_mg_{{jn}}-
    K^i_ng_{{jm}}\big)\\&
    -\frac1{N-2}\big(\delta^i_m{\overline K}{}_{jn}-\delta^i_n
    {\overline K}{}_{jm}+
    {\overline K}{}^i_m\overline g_{{jn}}-
    {\overline K}{}^i_n\overline g_{{jm}}
    \big)\\&-
    \frac{u}{N-2}\big(\delta^i_mT^\alpha_{{jn};\alpha}-
    \delta^i_nT^\alpha_{{jm};\alpha}-
    \delta^i_m\overline
    T^\alpha_{{jn}|\alpha}+\delta^i_n\overline
    T^\alpha_{{jm}|\alpha}\big)\\&-
    \frac{u}{N-2}\big(g^{{i\alpha}}
    T^\beta_{{\alpha m};\beta}g_{{jn}}-
    g^{{i\alpha}}
    T^\beta_{{\alpha n};\beta}g_{{jm}}
    -\overline g^{{i\alpha}}\overline T^\beta_{{\alpha m}
    |\beta}\overline g_{{jn}}+
    \overline g^{{i\alpha}}
    \overline T^\beta_{{\alpha n}|\beta}
    \overline g_{{jm}}
    \big)\\&-
    \frac{v'+w}{N-2}\big(\delta^i_mT^\alpha_{{j\beta}}
    T^\beta_{{\alpha n}}
    -\delta^i_nT^\alpha_{{j\beta}}T^\beta_{{\alpha
    m}}-\delta^i_m\overline T^\alpha_{{j\beta}}\overline
    T^\beta_{{\alpha n}}+\delta^i_n\overline
    T^\alpha_{{j\beta}}\overline T^\beta_{{\alpha
    m}}
    \big)\\&-
    \frac{v'+w}{N-2}\big(g^{{i\alpha}}
    T^\beta_{{\alpha\gamma}}T^\gamma_{{\beta
    m}}g_{{jn}}-g^{{i\alpha}}
    T^\beta_{{\alpha\gamma}}T^\gamma_{{\beta
    n}}g_{{jm}}-\overline g^{{i\alpha}}\overline
    T^\beta_{{\alpha\gamma}}
    \overline T^\gamma_{{\beta m}}\overline
    g_{{jn}}+\overline g^{{i\alpha}}\overline
    T^\beta_{{\alpha\gamma}}\overline T^\gamma_{{\beta n}}
    \overline
    g_{{jm}}
    \big)\\&+
    \frac1{(N-1)(N-2)}\Big(K(\delta^i_mg_{{jn}}-\delta^i_ng_{{jm}})
    -{\overline K}(\delta^i_m\overline g_{{jn}}-\delta^i_n\overline g_{{jm}})\Big)
    \\&-\frac{v'+w}{(N-1)(N-2)}\Big(T^\alpha_{{\gamma\beta}}
    T^\beta_{{\alpha\delta}}
    g^{{\gamma\delta}}\big(\delta^i_mg_{{jn}}
    -\delta^i_ng_{{jm}}\big)-
    \overline T^\alpha_{{\gamma\beta}}
    \overline T^\beta_{{\alpha\delta}}
    \overline g^{{\gamma\delta}}\big(\delta^i_m
    \overline g_{{jn}}-\delta^i_n\overline g_{{jm}}\big)\Big),
    \endaligned\label{eq:R0Rfaktorizacija}
  \end{equation}

  \noindent for covariant derivation with regard to affine
  connection of the associated space $\overline{\mathbb R}_N$
  denoted by vertical line $|$.

It is evident that geometrical objects $\tau^i_{(r)jk}$ and
$\overline\tau^i_{(r)jk}$ given by the equations (\ref{eq:ttauc},
\ref{eq:tauc}) are anti-symmetric by indices $j$ and $k$. After
antisymmetrize the invariants $\mathcal T^i_{(r)jk}$ by these
indices, we obtain that it holds

\begin{equation}
  \aligned
  (\overline T^i_{jm}&+\overline\tau^i_{(r)jm})_{|n}-
  (T^i_{jm}+\tau^i_{(r)jm})_{;n}=\overline T^i_{jm|n}-
  T^i_{jm;n}+\overline\tau^i_{(r)jm|n}-
  \tau^i_{(r)jm;n}\\&=
  \overline\zeta^i_{(s_1)\alpha n}(\overline T^\alpha_{jm}+
  \overline\tau^\alpha_{(r)jm})-
  \overline\zeta^\alpha_{(s_2)jn}(\overline T^i_{\alpha m}+
  \overline\tau^i_{(r)\alpha m})-
  \overline\zeta^\alpha_{(s_3)mn}(\overline T^i_{j\alpha}+
  \overline\tau^i_{(r)j\alpha})\\&-
  \zeta^i_{(s_1)\alpha n}(T^\alpha_{jm}+\tau^\alpha_{(r)jm})+
  \zeta^\alpha_{(s_2)jn}(T^i_{\alpha m}+\tau^i_{(r)\alpha m})+
  \zeta^\alpha_{(s_3)mn}(T^i_{j\alpha}+\tau^i_{(r)j\alpha}),
  \endaligned
\end{equation}

\noindent for $s=(s_1,s_2,s_3)\in\{1,2\}^3$ and the above defined
$\zeta^i_{(s_k)jm}$. Based on the last of previous equalities, we
proved that is

  \begin{align}
    &\overline T^i_{{jm}|n}=
    T^i_{{jm};n}+\overline\sigma{}^i_{(s)(r)jmn}-
    {\sigma}^i_{(s)(r)jmn},\label{eq:gammat1}
  \end{align}

  \noindent where is

  \begin{align*}
    \sigma{}^i_{(s)(r)jmn}=\tau^i_{(r)jm;n}-
  \zeta^i_{(s_1)\alpha n}(T^\alpha_{jm}+\tau^\alpha_{(r)jm})+
  \zeta^\alpha_{(s_2)jn}(T^i_{\alpha m}+\tau^i_{(r)\alpha m})+
  \zeta^\alpha_{(s_3)mn}(T^i_{j\alpha}+\tau^i_{(r)j\alpha}),
  \end{align*}

  \noindent and the corresponding
  ${\overline\sigma}^i_{(s)(r)jmn},s\in\{1,2\}^3,r\in\{1,2\}^5$.
  Furthermore, based on the equalities

  \begin{equation*}
    \aligned
    \big(\overline T^\alpha_{{jm}}+
    {\overline\tau}^\alpha_{(r^1)jm}\big)\big(
    \overline T^i_{{\alpha n}}+
    {\overline\tau}^i_{(r^2)\alpha n}\big)-
    \big(T^\alpha_{{jm}}+
    {\tau}^\alpha_{(r^1)jm}\big)\big(
    T^i_{{\alpha n}}+
    {\tau}^i_{(r^2)\alpha n}\big)=0,
    \endaligned
  \end{equation*}

  \noindent $r^1,r^2\in\{1,2\}^5$, we conclude that it is satisfied

  \begin{equation}
    \overline T^\alpha_{{jm}}
    \overline T^i_{{\alpha n}}=
    T^\alpha_{{jm}}
    T^i_{{\alpha n}}+
    \Theta^i_{(r^1)(r^2)jmn}-
    {\overline\Theta}^i_{(r^1)(r^2)jmn},
    \label{eq:ttTott}
  \end{equation}

  \noindent for

  \begin{align*}
    &\Theta^i_{(r^1)(r^2)jmn}=
    T^\alpha_{{jm}}\tau^i_{(r^2)\alpha
    n}+T^i_{{\alpha
    n}}\tau^\alpha_{(r^1)jm}+\tau^\alpha_{(r^1)jm}
    \tau^i_{(r^2)\alpha n},
  \end{align*}

  \noindent and the corresponding
  ${\overline\Theta}^i_{(r^1)(r^2)jmn}$.

  Based on the equations (\ref{eq:R}, \ref{eq:R0Rfaktorizacija} --
  \ref{eq:ttTott}), we obtain that it holds the
  relation

  \begin{equation}
    \aligned
    {\overline K}{}^i_{jmn}&=K^i_{jmn}+
    \frac1{N-2}\big(\delta^i_mK_{jn}-\delta^i_n
    K_{jm}+ K^i_mg_{{jn}}-
    K^i_ng_{{jm}}\big)\\&
    -\frac1{N-2}\big(\delta^i_m{\overline K}_{jn}-\delta^i_n
    {\overline K}_{jm}+
    {\overline K}^i_m\overline g_{{jn}}-
    {\overline K}{}^i_n\overline g_{{jm}}
    \big)\\&-
    \frac{u}{N-2}\big(\delta^i_mT^\alpha_{{jn};\alpha}-
    \delta^i_nT^\alpha_{{jm};\alpha}-
    \delta^i_m\overline
    T^\alpha_{{jn}|\alpha}+\delta^i_n\overline
    T^\alpha_{{jm}|\alpha}\big)\\&-
    \frac{u}{N-2}\big(g^{{i\alpha}}
    T^\beta_{{\alpha m};\beta}g_{{jn}}-
    g^{{i\alpha}}
    T^\beta_{{\alpha n};\beta}g_{{jm}}
    -\overline g^{{i\alpha}}\overline T^\beta_{{\alpha m}
    |\beta}\overline g_{{jn}}+
    \overline g^{{i\alpha}}
    \overline T^\beta_{{\alpha n}|\beta}
    \overline g_{{jm}}
    \big)\\&-
    \frac{v'+w}{N-2}\big(\delta^i_mT^\alpha_{{j\beta}}
    T^\beta_{{\alpha n}}
    -\delta^i_nT^\alpha_{{j\beta}}T^\beta_{{\alpha
    m}}-\delta^i_m\overline T^\alpha_{{j\beta}}\overline
    T^\beta_{{\alpha n}}+\delta^i_n\overline
    T^\alpha_{{j\beta}}\overline T^\beta_{{\alpha
    m}}
    \big)\\&-
    \frac{v'+w}{N-2}\big(g^{{i\alpha}}
    T^\beta_{{\alpha\gamma}}T^\gamma_{{\beta
    m}}g_{{jn}}-g^{{i\alpha}}
    T^\beta_{{\alpha\gamma}}T^\gamma_{{\beta
    n}}g_{{jm}}-\overline g^{{i\alpha}}\overline
    T^\beta_{{\alpha\gamma}}
    \overline T^\gamma_{{\beta m}}\overline
    g_{{jn}}+\overline g^{{i\alpha}}\overline
    T^\beta_{{\alpha\gamma}}\overline T^\gamma_{{\beta n}}
    \overline
    g_{{jm}}
    \big)\\&+
    \frac1{(N-1)(N-2)}\Big(K(\delta^i_mg_{{jn}}-\delta^i_ng_{{jm}})
    -{\overline K}(\delta^i_m\overline g_{{jn}}-\delta^i_n\overline g_{{jm}})\Big)
    \\&-\frac{v'+w}{(N-1)(N-2)}\Big(T^\alpha_{{\gamma\beta}}
    T^\beta_{{\alpha\delta}}
    g^{{\gamma\delta}}\big(\delta^i_mg_{{jn}}
    -\delta^i_ng_{{jm}}\big)-
    \overline T^\alpha_{{\gamma\beta}}
    \overline T^\beta_{{\alpha\delta}}
    \overline g^{{\gamma\delta}}\big(\delta^i_m
    \overline g_{{jn}}-\delta^i_n\overline g_{{jm}}\big)\Big)\\&
    -u\sigma^i_{(s^1)(r^1)jmn}-
    u'\sigma^i_{(s^2)(r^2)jnm}+v\Theta^i_{(r^3)(r^4)jmn}+
    v'\Theta^i_{(r^5)(r^6)jnm}+w\Theta^i_{(r^7)(r^8)mnj}\\&+
    u\overline\sigma^i_{(s^1)(r^1)jmn}+
    u'\overline\sigma^i_{(s^2)(r^2)jnm}-v\overline\Theta^i_{(r^3)(r^4)jmn}-
    v'\overline\Theta^i_{(r^5)(r^6)jnm}-w\overline\Theta^i_{(r^7)(r^8)mnj},
  \endaligned
  \end{equation}

  \noindent
  for $s^1,s^2\in\{1,2\}^3,r^1,\ldots,r^8\in\{1,2\}^5$. From this
  equation, we obtain that is

  \begin{equation*}
    {\overline C}{}^i_{(\rho)jmn}={C}^i_{(\rho)jmn},
  \end{equation*}

  \noindent for

  \begin{equation}
  \aligned
    C^i_{(\rho)jmn}&=K^i_{jmn}+
      \frac1{N-2}(\delta^i_mK_{jn}-\delta^i_nK_{jm}+
      K^i_mg_{{jn}}-K^i_ng_{{jm}})\\&+
      \frac{K}{(N-1)(N-2)}(\delta^i_mg_{{jn}}-
      \delta^i_ng_{{jm}})-u\sigma^i_{(s^1)(r^1)jmn}-
    u'\sigma^i_{(s^2)(r^2)jnm}\\&+v\Theta^i_{(r^3)(r^4)jmn}+
    v'\Theta^i_{(r^5)(r^6)jnm}+w\Theta^i_{(r^7)(r^8)mnj}\\&-
      \frac
      u{N-2}\big(\delta^i_mT^\alpha_{{jn};\alpha}-
      \delta^i_nT^\alpha_{{jm};\alpha}+
      g^{{i\alpha}}T^\beta_{{\alpha
      m};\beta}g_{{jn}}-
      g^{{i\alpha}}T^\beta_{{\alpha
      n};\beta}g_{{jm}}
      \big)\\&-
      \frac{v'+w}{N-2}\big(\delta^i_mT^\alpha_{{j\beta}}
      T^\beta_{{\alpha n}}-
      \delta^i_nT^\alpha_{{j\beta}}
      T^\beta_{{\alpha m}}+
      g^{{i\alpha}}T^\beta_{{\alpha\gamma}}
      T^\gamma_{{\beta m}}g_{{jn}}-
      g^{{i\alpha}}T^\beta_{{\alpha\gamma}}
      T^\gamma_{{\beta n}}g_{{jm}}
      \big)\\&-
      \frac{v'+w}{(N-1)(N-2)}T^\alpha_{{\gamma\beta}}
      T^\beta_{{\alpha\delta}}g^{{\gamma\delta}}
      (\delta^i_mg_{{jn}}-\delta^i_ng_{{jm}})
  \endaligned\label{eq:cWeylinv}
  \end{equation}

  \noindent and $\rho=(s^1,s^2,\rho^1,\ldots,\rho^8)$.

  \pagebreak

  It holds the following theorem:

  \begin{thm}
    Let $f:\mathbb{GR}_N\rightarrow\mathbb{G\overline R}_N$ be a
    conformal mapping of generalized Riemannian space
    $\mathbb{GR}_N$. The families $C^i_{(\rho)jmn}$
    from the equation \emph{(\ref{eq:cWeylinv})} are families of
     invariants of the mapping $f$.\qed
  \end{thm}

  \begin{cor}
    Let $f:\mathbb{GR}_N\rightarrow\mathbb{G\overline R}_N$ be a
    conformal mapping of generalized Riemannian space
    $\mathbb{GR}_N$. Geometrical objects

    \begin{eqnarray}
  \aligned
    C_{(\rho)ijmn}&=K_{ijmn}+
      \frac1{N-2}(K_{jn}g_{im}-K_{jm}g_{in}+
      K_{im}g_{{jn}}-K_{in}g_{{jm}})\\&+
      \frac{K}{(N-1)(N-2)}(g_{im}g_{{jn}}-
      g_{in}g_{{jm}})-u\sigma_{(s^1)(r^1)ijmn}-
    u'\sigma_{(s^2)(r^2)ijnm}\\&+v\Theta_{(r^3)(r^4)ijmn}+
    v'\Theta_{(r^5)(r^6)ijnm}+w\Theta_{(r^7)(r^8)imnj}\\&-
      \frac
      u{N-2}\big(T^\alpha_{{jn};\alpha}g_{im}-
      T^\alpha_{{jm};\alpha}g_{in}+
      T^\beta_{{im};\beta}g_{{jn}}-
      T^\beta_{{in};\beta}g_{{jm}}
      \big)\\&-
      \frac{v'+w}{N-2}\big(T^\alpha_{{j\beta}}
      T^\beta_{{\alpha n}}g_{im}-
      T^\alpha_{{j\beta}}
      T^\beta_{{\alpha m}}g_{in}+
      T^\beta_{{i\gamma}}
      T^\gamma_{{\beta m}}g_{{jn}}-
      T^\beta_{{i\gamma}}
      T^\gamma_{{\beta n}}g_{{jm}}
      \big)\\&-
      \frac{v'+w}{(N-1)(N-2)}T^\alpha_{{\gamma\beta}}
      T^\beta_{{\alpha\delta}}g^{{\gamma\delta}}
      (g_{im}g_{{jn}}-g_{in}g_{{jm}}),
  \endaligned\label{eq:cconfWeylinv}
  \end{eqnarray}

  \noindent for
  $\sigma_{(s)(r)ijmn}=g_{i\alpha}\sigma^\alpha_{(s)(r)jmn}$ and
  $\Theta_{(r^i)(r^j)ijmn}=g_{i\alpha}\Theta^\alpha_{(r^i)(r^j)jmn}$,
  are invariants of the mapping $f$.\qed
  \end{cor}

  \noindent The above obtained
  invariants $C^i_{(\rho)jmn}$ and
  $C_{(\rho)ijmn}$ are \textbf{the $\rho$-th
  conformal} and \textbf{the $\rho$-th
  conformal covariant  invariants of Weyl type} of the conformal mapping $f$.

  \section{Conclusion}

  In this paper, we studied transformations of Christoffel symbols
  and curvature tensors of a generalized Riemannian space under conformal mappings.
  From these transformations, we obtained invariants of these
  transformations. Furthermore, invariance of
  some of the invariants $\mathcal T^i_{(r)jk}$ from the equation
  (\ref{eq:Thomas}) is necessary and sufficient condition for a
  mapping\linebreak $f:\mathbb{GR}_N\rightarrow\mathbb{G\overline R}_N$ to be
  a conformal one.

  In further research, motivated by the results from this paper, we will obtain
  invariants of $F$-planar mappings of non-symmetric affine
  connection spaces.

  \section*{Acknowledgements}

  This paper is financially supported by Serbian Ministry of
  Education, Science and Technological Development, Grant No.
  174012.

  \pagebreak

%\noindent\textbf{Author:}

%\texttt{Nenad O. Vesi\'c}

%\texttt{Faculty of Science and Mathematics}

%adress: \texttt{Vi\v segradska 33, 18000 Ni\v s, Serbia}

%Contact e-mail: \texttt{n.o.vesic@outlook.com}

\end{document}